# A NON-TRIVIAL $(\mathbf{R},+)$ PRINCIPAL BUNDLE OVER A CONTRACTIBLE BASE?

PATRICK IGLESIAS-ZEMMOUR

ABSTRACT. We investigate the properties of a specific quotient space construction, the "warped projection" $\pi : \mathrm{W}_\alpha \to \mathrm{D}_\alpha$, over a smoothly contractible base. In a previous version of this work, it was claimed that this structure constituted a non-trivial principal bundle. We revisit this claim and observe that, due to the existence of flat functions in the smooth category, the projection fails to satisfy the strict condition of local triviality along the plots, required for diffeological bundles. However, the structure remains rich: it possesses a smooth, free, and fiber-transitive group action. Drawing on the concept of *vector pseudo-bundles* introduced by Christensen and Wu, we propose that this object is best understood as a non-trivial **principal pseudo-bundle**. This example thus serves to clarify the boundary between strict bundles and generalized pseudo-bundles in the context of singular base spaces.

## INTRODUCTION

The theory of diffeology was developed to extend the methods of differential geometry to spaces that are not manifolds, such as quotients, singularities, and infinite-dimensional groups. A central question in this framework is the classification of fiber bundles. While the theory of *diffeological fiber bundles* is well-established [2, 3], recent developments have suggested that for certain singular spaces, the strict axioms of local triviality may be too restrictive to capture all geometrically relevant structures.

In the context of tangent bundles, Christensen and Wu [1] and Wu [4] have introduced the notion of **diffeological vector pseudo-bundles**. These are smooth surjections with vector space fibers and smooth operations, but which do not necessarily satisfy the condition of local triviality along plots. They arise naturally when considering the tangent spaces of singular objects.

In this paper, we investigate the non-abelian (or principal) analogue of this phenomenon. We construct a specific object, $\pi : \mathrm{W}_\alpha \to \mathrm{D}_\alpha$, which we call the **warped projection**. The base space $\mathrm{D}_\alpha$ is the quotient of the disk by an irrational rotation, a space known to be smoothly contractible. The total space $\mathrm{W}_\alpha$ is constructed by "warping" the product with $\mathbf{R}$ via a specific action.

**Correction and Reinterpretation.** In a preliminary version of this manuscript, it was asserted that $\pi : \mathrm{W}_\alpha \to \mathrm{D}_\alpha$ was a principal bundle. As we demonstrate in Section 2, this claim was incorrect. The subtle existence of flat functions (smooth functions with vanishing derivatives of all orders) prevents the inverse of the action map from being smooth. Consequently, the structure fails the *induction criterion* required for principal bundles.

However, this failure is instructive. The warped projection satisfies all the conditions of a bundle *except* the smoothness of the division map. It carries a smooth, free, and transitive action of $(\mathbf{R},+)$ on the fibers. It is, therefore, the principal equivalent of a vector pseudo-bundle. Furthermore, we show that this pseudo-bundle is **non-trivial**: it is not isomorphic to the product $\mathrm{D}_\alpha \times \mathbf{R}$ in the category of pseudo-bundles.

This result highlights a divergence between the topological and diffeological categories. While a principal bundle over a contractible base must be trivial, a principal *pseudo-bundle* need not be. The warped projection $\mathrm{W}_\alpha$ thus stands as a foundational counter-example, delimiting the scope of bundle theory and motivating the further study of pseudo-bundles over singular spaces.

## 1. The Construction

Let $\mathrm{D} = \{z \in \mathbf{C} \mid |z| \leq 1\}$ be the closed unit disk equipped with the subset diffeology of $\mathbf{C}$. Let $\alpha \in \mathbf{R} \setminus \mathbf{Q}$ be a fixed irrational number. We consider the action of $\mathbf{Z}$ on D by rotation: $n \cdot z = z e^{i2\pi n\alpha}$. The quotient $\mathrm{D}_\alpha = \mathrm{D}/\mathbf{Z}$ is the base space.

We lift this action to the product $\mathrm{D} \times \mathbf{R}$ via the formula:

$$n \cdot (z, t) = (z e^{i2\pi n\alpha}, t + n|z|^2). \tag{1}$$

This action is smooth. Let $\mathrm{W}_\alpha = (\mathrm{D} \times \mathbf{R})/\mathbf{Z}$ be the quotient space equipped with the quotient diffeology. The projection $\mathrm{pr}_1 : \mathrm{D} \times \mathbf{R} \to \mathrm{D}$ is equivariant and descends to a smooth map $\pi : \mathrm{W}_\alpha \to \mathrm{D}_\alpha$. We call this map the **warped projection**.

Furthermore, there is a natural action of $(\mathbf{R},+)$ on $\mathrm{W}_\alpha$ defined by $s \cdot [z, t] = [z, t+s]$. This action is smooth, free, and transitive on the fibers of $\pi$.

## 2. The Smooth Obstruction

Recall that for $\pi$ to be a principal bundle, the map $\mathrm{F} : \mathrm{W}_\alpha \times \mathbf{R} \to \mathrm{W}_\alpha \times_{\mathrm{D}_\alpha} \mathrm{W}_\alpha$ defined by $\mathrm{F}(y, s) = (y, s \cdot y)$ must be an induction (a diffeomorphism onto its image). While F is injective, we show here that its inverse is **not smooth** in the standard $\mathscr{C}^\infty$ diffeology.

**Proposition 1.** *The map $\mathrm{F}^{-1}$ is not smooth. Consequently, $\pi : \mathrm{W}_\alpha \to \mathrm{D}_\alpha$ is not a diffeological principal bundle.*

*Proof.* We construct a smooth plot in $\mathrm{W}_\alpha \times_{\mathrm{D}_\alpha} \mathrm{W}_\alpha$ whose lift to $\mathbf{R}$ via $\mathrm{F}^{-1}$ is discontinuous. This is the **rapidly oscillating plot (ROP)**.

Let $w_1 : \mathbf{R} \to \mathrm{D}$ be a smooth function such that $w_1(t)$ is real, non-negative, and possesses a sequence of zeros $t_k \to 0$ ($t_k \neq 0$). Assume $w_1$ is *flat* at every $t_k$ (all derivatives vanish).

Let $\mathrm{I}_k$ be an open interval around $t_k$ such that the $\mathrm{I}_k$ are disjoint. Define a sequence of integers $\mathrm{N}_k$ such that $\mathrm{N}_k \to \infty$ as $k \to \infty$. Specifically, we can choose $\mathrm{N}_k$ large enough such that $\mathrm{N}_k \sup_{t \in \mathrm{I}_k} |w_1(t)|^2 \geq 1$.

Define a function $m : \mathbf{R} \to \mathbf{Z}$ such that $m(t) = \mathrm{N}_k$ if $t \in \mathrm{I}_k$, and $m(t) = 0$ otherwise. Define the function $s : \mathbf{R} \to \mathbf{R}$ by $s(t) = m(t)|w_1(t)|^2$. As shown previously, $s(t)$ is discontinuous at $t = 0$ (oscillating between 0 and $\geq 1$), so it is not smooth.

Consider the path in the domain $\mathrm{W}_\alpha \times \mathbf{R}$:

$$\Gamma(t) = ([w_1(t), 0], s(t)).$$

Since the second component $s(t)$ is not smooth, $\Gamma$ is not a plot of $\mathrm{W}_\alpha \times \mathbf{R}$.

Now consider the image path $\mathrm{P}(t) = \mathrm{F}(\Gamma(t))$ in the codomain. Recall that $\mathrm{F}(y, s) = (y, s \cdot y)$. Here $y = [w_1(t), 0]$. The action is $s \cdot [z, \tau] = [z, \tau + s]$. Thus:

$$\mathrm{P}(t) = ([w_1(t), 0], [w_1(t), s(t)]).$$

We must show that $\mathrm{P}(t)$ is a smooth plot of $\mathrm{W}_\alpha \times \mathrm{W}_\alpha$. The first component $t \mapsto [w_1(t), 0]$ is clearly smooth. For the second component, we claim that the map $t \mapsto [w_1(t), s(t)]$ is smooth.

Consider the lift $\tilde{\gamma} : \mathbf{R} \to \mathrm{D} \times \mathbf{R}$ defined by:

$$\tilde{\gamma}(t) = (w_1(t) e^{-i2\pi m(t)\alpha}, 0).$$

This map is smooth because $w_1$ is flat at the points where $m(t)$ is discontinuous; the flatness suppresses the phase jumps. Now, let us check its projection. Using the equivalence relation $(z, \tau) \sim (n \cdot z, \tau + n|z|^2)$ with $n = m(t)$:

$$\begin{aligned} m(t) \cdot \tilde{\gamma}(t) &= m(t) \cdot (w_1(t) e^{-i2\pi m(t)\alpha}, 0) \\ &= (w_1(t) e^{-i2\pi m(t)\alpha} e^{i2\pi m(t)\alpha}, 0 + m(t)|w_1(t) e^{-i2\pi m(t)\alpha}|^2) \\ &= (w_1(t), m(t)|w_1(t)|^2) \\ &= (w_1(t), s(t)). \end{aligned}$$

Therefore, $\varpi(\tilde{\gamma}(t)) = [w_1(t), s(t)]$. Since $[w_1(t), s(t)]$ admits a smooth lift $\tilde{\gamma}$, it is a smooth plot in $\mathrm{W}_\alpha$.

**Conclusion:** F maps a non-smooth path $\Gamma$ to a smooth path P. Therefore, F does not reflect smoothness; it is not an induction. □

## 3. The Principal Pseudo-Bundle Structure

Following the terminology of [1], we define a **principal pseudo-bundle** as a smooth surjection $\pi : \mathrm{Y} \to \mathrm{X}$ equipped with a smooth, free, and fiber-transitive action of a diffeological group G. Crucially, we do not require the action map to be an induction, only a smooth map.

**Proposition 2.** *The warped projection* $\pi : \mathrm{W}_\alpha \to \mathrm{D}_\alpha$ *is a principal* $(\mathbf{R}, +)$ *pseudo-bundle.*

*Proof.* The projection $\pi$ is smooth by construction. The action of $\mathbf{R}$ on $\mathrm{W}_\alpha$ is smooth because it descends from the smooth action on $\mathrm{D} \times \mathbf{R}$. It is free and transitive on the fibers (as shown in the injectivity of F). Thus, it satisfies the axioms of a pseudo-bundle. □

**Theorem 3** (Non-Triviality)**.** *The pseudo-bundle* $\pi : \mathrm{W}_\alpha \to \mathrm{D}_\alpha$ *is not isomorphic to the trivial pseudo-bundle* $\mathrm{pr}_1 : \mathrm{D}_\alpha \times \mathbf{R} \to \mathrm{D}_\alpha$.

*Proof.* An isomorphism of pseudo-bundles must be a diffeomorphism. If $\Phi : D_\alpha \times \mathbf{R} \to W_\alpha$ were an isomorphism, its inverse $\Phi^{-1}$ would be smooth. The map $\Phi$ corresponds to the action map F restricted to a global section (which exists in the trivial case). However, we have proven that $F^{-1}$ is not smooth. Specifically, the "division" operation required to map an element of $W_\alpha$ back to the trivial product coordinates involves the same discontinuous function $s(t)$ constructed in the ROP. Since the inverse map is not smooth, the spaces are not diffeomorphic, and the pseudo-bundle is non-trivial. □

## Concluding Remarks

The warped projection $W_\alpha$ serves as a boundary marker in the theory of diffeological bundles. In the category of **Analytic Diffeology**, where flat functions do not exist, the ROP obstruction vanishes, and $W_\alpha$ becomes a true principal bundle. In the standard **Smooth Diffeology**, the ROP obstruction prevents it from being a bundle, demoting it to a pseudo-bundle. This distinction is invisible to topology but fundamental to diffeology, illustrating the theory's capacity to detect fine analytic structures.

Einstein Institute of Mathematics, The Hebrew University of Jerusalem, Israel
*Email address*: piz@math.huji.ac.il